\documentclass[11pt]{article}
\usepackage{latexsym}
\usepackage{amsfonts}
\usepackage{mathrsfs,amsthm,amsmath}
\usepackage{color}
\usepackage{todonotes}
\newtheorem{theorem}{Theorem}[section]
\newtheorem{proposition}[theorem]{Proposition}

\newtheorem{remark}{Remark}[section]
\newtheorem{example}{Example}[section]
\newtheorem{condition}{Condition}[section]
\begin{document}
\title{Some vector-valued examples of noncentral moderate deviation 
results\thanks{The authors acknowledge the support of MUR Excellence Department Project awarded to the 
Department of Mathematics, University of Rome Tor Vergata (CUP E83C23000330006), of University of Rome 
Tor Vergata (CUP E83C25000630005) Research Project METRO, and of INdAM-GNAMPA.}}
\author{Claudio Macci\thanks{Address: Dipartimento di Matematica,
Università di Roma Tor Vergata, Via della Ricerca Scientifica,
I-00133 Rome, Italy. e-mail: \texttt{macci@mat.uniroma2.it}} \and
Barbara Pacchiarotti\thanks{Address: Dipartimento di Matematica,
Università di Roma Tor Vergata, Via della Ricerca Scientifica,
I-00133 Rome, Italy. e-mail: \texttt{pacchiar@mat.uniroma2.it}}}
\maketitle
\begin{abstract}
	The term \emph{noncentral moderate deviations} is used in the literature to mean a class of large 
	deviation principles that, in some sense, fills the gap between the convergence in probability to a 
	constant (governed by a reference large deviation principle) and a weak convergence to a non-Gaussian
	(and non-degenerating) distribution. Several examples can be found in the literature, mainly for 
	real-valued random variables (see, e.g.,~\cite{GiulianoMacci} and the references cited therein). In this 
	paper we present some examples with vector-valued random variables.
\end{abstract}

\section{Introduction}
The term \emph{moderate deviations} is used when we have a suitable class of large deviation principles which fills the gap
between a convergence to a constant, governed by a \emph{reference large deviation principle}, and an asymptotic normality
result (see, e.g.,~\cite{DemboZeitouni} as a reference on large deviations). 
A very well-known example is provided by Theorem 3.7.1 in \cite{DemboZeitouni} which provides a class of large deviation 
principles which fills the gap between the convergence to a constant (provided by the \emph{Law of the Large Numbers}) and
the weak convergence to a Gaussian distribution (provided by the \emph{Central Limit Theorem}). More recently the term 
\emph{noncentral moderate deviations} has been used when one has a class of large deviation principles of the same kind, but 
the weak convergence is towards a non-Gaussian (and non-degenerating) distribution. Several examples of noncentral moderate 
deviation results with real-valued random variables can be found in \cite{GiulianoMacci} (see also the references cited therein).
In that reference the examples concern cases in which one has the weak convergence towards Gumbel, exponential and
Laplace distributions.

The aim of this paper is to present some \emph{vector-valued} examples of noncentral moderate deviation results. More
precisely we mean examples with $\mathbb{R}^h$-valued random variables. Some examples with time-changed multivariate Lévy processes
with linear combinations of inverse stable subordinators appear in \cite{GuptaMacci}. Another example appears in 
\cite{LeonenkoMacciPacchiarotti} (Section 3) but, in that case, the weak convergence is trivial (indeed the involved 
random variables are identically distributed). We can also say that, in general, we can obtain univariate versions of the examples
in this paper by setting $h=1$ (actually one should consider the case $h=2$ in Example \ref{ex2} because that example
concerns $\mathbb{R}^{h-1}$-valued random variables). The noncentral moderate deviation results will be given in Propositions 
\ref{prop:ncMD-IMM}, \ref{prop:ncMD-Levy}, \ref{prop:ncMD-Poisson} and \ref{prop:ncMD-continuity} (see also Remarks 
\ref{rem:fills-the-gap-IMM}, \ref{rem:fills-the-gap-Levy}, \ref{rem:fills-the-gap-Poisson} and 
\ref{rem:fills-the-gap-contraction-continuity}).
We remark that some weak convergence results in this paper (i.e. Propositions \ref{prop:weak-convergence-IMM}, 
\ref{prop:weak-convergence-Levy} and \ref{prop:weak-convergence-Poisson}) are proved by checking the pointwise convergence of the 
moment generating functions (and not of characteristic functions); for instance one could refer to Theorem 2 in 
\cite{MukherjeaRaoSuen}.

In Section \ref{sec:preliminaries} we recall some preliminaries. Section \ref{sec:IMM} is devoted to multivariate Lévy processes 
time-changed with an inverse stable subordinator. In this way we generalize the results for univariate Lévy processes presented 
in \cite{IulianoMacciMeoli} (Section 3). In Section \ref{sec:Levy} we consider multivariate Lévy processes time-changed with a 
light-tailed independent subordinator. In Section \ref{sec:Poisson} we present the results for a class of examples in which the
random variables converge weakly to a compound Poisson distributed law. In particular, as explained in Remark 
\ref{rem:on-the-well-known-case-Poisson}, this example provides a generalization of the well-known weak convergence of 
Binomial distributions to the Poisson distribution. Interestingly, in Sections \ref{sec:Levy} and \ref{sec:Poisson}, we have two 
different inequalities for the rate functions $I_{\mathrm{LD}}$ and $I_{\mathrm{MD}}$ (see eqs.~\eqref{eq:inequality-Levy} and \eqref{eq:inequality-Poisson} below).

In the final Section \ref{sec:contraction-continuity} we present the results for a class of examples which is obtained by
considering the application of continuous mappings to the random variables in Theorem 3.7.1 in \cite{DemboZeitouni}
recalled above. Indeed, we can do that by considering suitable applications of the Contraction Principle (see, e.g.,~Theorem 4.2.1 in 
\cite{DemboZeitouni}) and the Continuity Theorem for the weak convergence. Moreover, we remark that, obviously, this trick (i.e. the 
application of a continuous mapping to the random variables involved in a moderate deviation result) can be considered to obtain
other (possibly noncentral) moderate deviation results. In Section \ref{sec:contraction-continuity} we also study in detail two
examples. In Example \ref{ex1} we deal with random variables taking values on the simplex, and the weak convergence is towards the 
logistic Normal distribution (see, e.g.,~\cite{Aitchison}, Chapter 6). In Example \ref{ex2} we deal with a weak convergence towards 
a skew Normal distribution by considering a suitable continuous mapping in the literature (see, e.g.,~eqs.~(6) and (7) in 
\cite{AzzaliniCapitanioJRSSB} for the transformation method (b)). A reader interested to other examples could see 
\cite{MateufiguerasMontiEgozcue} for a wide source of examples of distributions on the simplex, and \cite{AzzaliniCapitanio}
as a reference on skew distributions.

\section{Preliminaries on large (and moderate) deviations}\label{sec:preliminaries}
We start with some basic definitions (see, e.g.,~\cite{DemboZeitouni}, pages 4-5). Let $\mathcal{X}$ be a topological 
space and let $\{B_n:n\geq 1\}$ be a sequence of $\mathcal{X}$-valued random variables defined on the same probability 
space $(\Omega,\mathcal{F},P)$.
A sequence $\{v_n:n\geq 1\}$ such that $v_n\to\infty$ (as $n\to\infty$) is called a \emph{speed function}, 
and a lower semicontinuous function $I:\mathcal{X}\to[0,\infty]$ is called a \emph{rate function}. Then the 
sequence $\{B_n:n\geq 1\}$ satisfies the large deviation principle (LDP from now on) with speed $v_n$ and
rate function $I$ if
$$\limsup_{n\to\infty}\frac{1}{v_n}\log P(B_n\in C)\leq-\inf_{x\in C}I(x)\quad\mbox{for all closed sets}\ C,$$
and
$$\liminf_{n\to\infty}\frac{1}{v_n}\log P(B_n\in O)\geq-\inf_{x\in O}I(x)\quad\mbox{for all open sets}\ O.$$
The rate function $I$ is said to be \emph{good} if, for every $\eta\geq 0$, the level set $\{x\in\mathcal{X}:I(x)\leq\eta\}$
is compact.

The main large deviation tool used in this paper is the \emph{G\"artner Ellis Theorem} (see, e.g.,~Theorem 2.3.6 in 
\cite{DemboZeitouni}; for our aim we restrict to the final statement $(c)$ which provides the full LDP). 
We denote by $\langle\cdot,\cdot\rangle$ the inner product in $\mathbb{R}^h$. Let $\{B_n:n\geq 1\}$ 
be a sequence of random variables such that, for all $\theta\in\mathbb{R}^h$, the limit
$$\Lambda(\theta):=\lim_{n\to\infty}\frac{1}{v_n}\log\mathbb{E}\left[e^{v_n\langle\theta,B_n\rangle}\right]$$
exists as an extended real number, and that $0\in(\mathcal{D}(\Lambda))^\circ$, where 
$(\mathcal{D}(\Lambda))^\circ$ is the interior of $\mathcal{D}(\Lambda):=\{\theta\in\mathbb{R}^h:\Lambda(\theta)<\infty\}$.
Then, if $\Lambda$ is essentially smooth and lower semicontinuous, the sequence $\{B_n:n\geq 1\}$ satisfies the LDP with 
speed $v_n$ and good rate function $I=\Lambda^*$, i.e. the Legendre transform of $\Lambda$. For completeness we recall that 
$\Lambda$ is \emph{essentially smooth} (see, e.g.,~Definition 2.3.5 in \cite{DemboZeitouni}) if $(\mathcal{D}(\Lambda))^\circ$
is non-empty, the function $\Lambda$ is differentiable throughout $(\mathcal{D}(\Lambda))^\circ$, and the function $\Lambda$
is \emph{steep} (i.e. $\|\nabla\Lambda(\theta_n)\|$ tends to infinity for every sequence 
$\{\theta_n:n\geq 1\}\subset(\mathcal{D}(\Lambda))^\circ$ that converges to a boundary point of $\mathcal{D}(\Lambda)$).

Finally, we recall the \emph{Contraction Principle} (see, e.g.,~Theorem 4.2.1 in \cite{DemboZeitouni}). Let $f:\mathcal{X}\to
\mathcal{Y}$ be a continuous function, where $\mathcal{X}$ and $\mathcal{Y}$ are two topological spaces. Assume that 
$\{B_n:n\geq 1\}$ is a sequence of $\mathcal{X}$-valued random variables which satisfies the LDP with speed $v_n$ and good rate 
function $I$. Then the sequence $\{f(B_n):n\geq 1\}$ satisfies the LDP with good rate function $J$ defined by
$$J(y):=\inf\{I(x):x\in\mathcal{X},\ f(x)=y\}.$$

\begin{remark}\label{rem:preliminaries}
	The preliminaries presented above concern sequences of random variables, indeed we have a discrete parameter $n$.
	On the contrary, for the results in Sections \ref{sec:IMM} and \ref{sec:Levy}, we refer to a version of these
	preliminaries (except the Cramér Theorem) with continuous parameter $t$, where $t$ tends to infinity.
\end{remark}

\section{Multivariate Lévy processes time-changed with an inverse stable subordinator}\label{sec:IMM}
The aim of this section is to present a multivariate version of the results in \cite{IulianoMacciMeoli} 
(Section 3) for the univariate case. We recall some preliminaries presented in that reference (see also
the references cited therein). Let $\{S(t):t\geq 0\}$ be a $\mathbb{R}^h$-valued Lévy 
process as in the following Condition \ref{cond:IMM-Levy} (in \cite{IulianoMacciMeoli} we have $h=1$).

\begin{condition}\label{cond:IMM-Levy}
	Let $\{S(t):t\geq 0\}$ be a $\mathbb{R}^h$-valued Lévy process such that 
	$\kappa_S(\theta):=\log\mathbb{E}\left[e^{\langle\theta,S(1)\rangle}\right]$ is finite when $\theta$ belongs to a 
	neighborhood of the origin $0\in\mathbb{R}^h$, where $\langle\cdot,\cdot\rangle$ is the inner product in 
	$\mathbb{R}^h$.
\end{condition}
In particular, if Condition \ref{cond:IMM-Levy} holds, it is well-known that the component-wise 
expectation of the random vector $S(1)=(S_1(1),\ldots,S_h(1))$, i.e. $(\mathbb{E}[S_1(1)],\ldots,\mathbb{E}[S_h(1)])$, 
coincides with $\nabla\kappa_S(0)$.

Moreover, for $\nu\in(0,1)$, let $\{L_\nu(t):t\geq 0\}$ be an inverse stable subordinator. It is well-known that
	$$\mathbb{E}[e^{\theta L_\nu(t)}]=E_\nu(\theta t^\nu)\ \mbox{for all}\ \theta\in\mathbb{R},$$
where 
$$E_\nu(x):=\sum_{k=0}^\infty\frac{x^k}{\Gamma(\nu k+1)}$$
is the Mittag-Leffler function (see, e.g.,~\cite{GorenfloKilbasMainardiRogosin}, eq. (3.1.1)). We always assume 
that $\{S(t):t\geq 0\}$ and $\{L_\nu(t):t\geq 0\}$ are independent. 

We recall some other preliminaries. In view of the applications of the G\"artner Ellis Theorem in the proofs 
of Propositions \ref{prop:reference-LDP-IMM} and \ref{prop:ncMD-IMM}, we have
\begin{equation}\label{eq:ML-asymptotics}
	\left\{\begin{array}{l}
		E_\nu(x)\sim\frac{e^{x^{1/\nu}}}{\nu}\ \mbox{as}\ x\to\infty\\
		\mbox{for $y<0$, we have}\ \frac{1}{x}\log E_\nu(yx)\to 0\ \mbox{as}\ x\to\infty
	\end{array}\right.
\end{equation}
(see, e.g.,~Proposition 3.6 in \cite{GorenfloKilbasMainardiRogosin} for the case $\alpha\in(0,2)$; indeed 
$\alpha$ in that reference coincides with $\nu$ in this paper). Moreover, we consider the functions
	$$\alpha(\nu):=\left\{\begin{array}{ll}
		1-\nu/2&\ \mbox{if}\ \nabla\kappa_S(0)=0\\
		1-\nu&\ \mbox{if}\ \nabla\kappa_S(0)\neq 0,
	\end{array}\right.$$
and
$$f_\nu(y):=\left\{\begin{array}{ll}
	y^{1/\nu}&\ \mbox{if}\ y\geq 0\\
	0&\ \mbox{otherwise}.
\end{array}\right.$$
Finally, let $Q$ be the Hessian matrix of $\kappa_S$ the origin $0\in\mathbb{R}^h$. It is well-known that,
if we set $S(t)=(S_1(t),\ldots,S_h(t))$, then
$$Q=(\mathrm{Cov}(S_i(1),S_j(1)))_{i,j\in\{1,\ldots,h\}}.$$
This matrix plays the role of $q=\kappa_S^{\prime\prime}(0)$ for the case $h=1$ in \cite{IulianoMacciMeoli}
(see Condition 1.1 in that reference).

We have the following results.

\begin{proposition}[Reference LDP]\label{prop:reference-LDP-IMM}
	Assume that $f_\nu\circ\kappa_S$ is an essentially smooth function. Then 
	$\left\{\frac{S(L_\nu(t))}{t}:t>0\right\}$ satisfies the LDP with speed $t$ and good rate function $I_{\mathrm{LD}}$ 
	defined by
	$$I_{\mathrm{LD}}(x):=\sup_{\theta\in\mathbb{R}^h}\{\langle\theta,x\rangle-f_\nu(\kappa_S(\theta))\}.$$
\end{proposition}
\begin{proof}
	We consider a straightforward application of the G\"artner Ellis Theorem, and it is similar to the proof of 
	Proposition 3.1 in \cite{IulianoMacciMeoli} concerning the case $h=1$. Indeed (in the equalities below we have
	infinity when $\kappa_S(\theta)$ is equal to infinity) we have
	$$\lim_{t\to\infty}\frac{1}{t}\log\mathbb{E}\left[e^{t\langle\theta,S(L_\nu(t))/t\rangle}\right]
	=\lim_{t\to\infty}\frac{1}{t}\log E_\nu(\kappa_S(\theta)t^\nu)=f_\nu(\kappa_S(\theta))
	\quad(\mbox{for all}\ \theta\in\mathbb{R}^h).$$
\end{proof}

\begin{proposition}[Weak convergence]\label{prop:weak-convergence-IMM}
	We have the following statements.
	\begin{itemize}
		\item If $\nabla\kappa_S(0)=0$, then $\{t^{\alpha(\nu)}\frac{S(L_\nu(t))}{t}:t>0\}$ converges weakly 
		to a random vector $Z$, say, such that, given $L_\nu(1)$, it is Gaussian distributed with mean 
		$0\in\mathbb{R}^h$ and covariance matrix $QL_\nu(1)$.
		\item If $\nabla\kappa_S(0)\neq 0$, then $\{t^{\alpha(\nu)}\frac{S(L_\nu(t))}{t}:t>0\}$ converges weakly to
		$\nabla\kappa_S(0)L_\nu(1)$.
	\end{itemize}
\end{proposition}
\begin{proof}
	The proof is similar to the proof of Proposition 3.2 in \cite{IulianoMacciMeoli} concerning the case $h=1$.
	In both cases $\nabla\kappa_S(0)=0$ and $\nabla\kappa_S(0)\neq 0$ we study suitable limits (as $t\to\infty$) in 
	terms of the moment generating functions.
	
	If $\nabla\kappa_S(0)=0$, then we have
	\begin{multline*}
		\mathbb{E}\left[e^{\langle\theta,t^{\alpha(\nu)}\frac{S(L_\nu(t))}{t}\rangle}\right]
		=\mathbb{E}\left[e^{\langle\theta,\frac{S(L_\nu(t))}{t^{\nu/2}}\rangle}\right]
		=E_\nu\left(\kappa_S\left(\frac{\theta}{t^{\nu/2}}\right)t^\nu\right)\\ 
		=E_\nu\left(\left(\frac{\langle\theta,Q\theta\rangle}{2t^\nu}+o\left(\frac{1}{t^\nu}\right)\right)t^\nu\right)
		\to E_\nu\left(\frac{\langle\theta,Q\theta\rangle}{2}\right)\ \mbox{for all}\ \theta\in\mathbb{R}^h.
	\end{multline*}
	Thus the desired weak convergence is proved noting that the moment generating function of $Z$ is
	$$\mathbb{E}\left[e^{\langle\theta,Z\rangle}\right]
	=\mathbb{E}\left[e^{\frac{\langle\theta,Q\theta\rangle}{2}L_\nu(1)}\right]
	=E_\nu\left(\frac{\langle\theta,Q\theta\rangle}{2}\right)\ \mbox{for all}\ \theta\in\mathbb{R}^h.$$
	
	If $\nabla\kappa_S(0)\neq 0$, then we have
	\begin{multline*}
		\mathbb{E}\left[e^{\langle\theta,t^{\alpha(\nu)}\frac{S(L_\nu(t))}{t}\rangle}\right]
		=\mathbb{E}\left[e^{\langle\theta,\frac{S(L_\nu(t))}{t^\nu}\rangle}\right]
		=E_\nu\left(\kappa_S\left(\frac{\theta}{t^\nu}\right)t^\nu\right)\\
		=E_\nu\left(\left(\frac{\langle\theta,\nabla\kappa_S(0)\rangle}{t^\nu}+o\left(\frac{1}{t^\nu}\right)\right)t^\nu\right)
		\to E_\nu(\langle\theta,\nabla\kappa_S(0)\rangle)\ \mbox{for all}\ \theta\in\mathbb{R}^h.
	\end{multline*}
	Thus the desired weak convergence is proved noting that
	$$\mathbb{E}\left[e^{\langle\theta,\nabla\kappa_S(0)L_\nu(1)\rangle}\right]
	=\mathbb{E}\left[e^{\langle\theta,\kappa_S(0)L_\nu(1)\rangle}\right]
	=E_\nu(\langle\theta,\nabla\kappa_S(0)\rangle)\ \mbox{for all}\ \theta\in\mathbb{R}^h.$$
\end{proof}

\begin{proposition}[Noncentral moderate deviations]\label{prop:ncMD-IMM}
	Assume that $Q$ is not the null matrix if $\nabla\kappa_S(0)=0$. 
	Then, for every family of positive numbers $\{a_t:t>0\}$ such that $a_t\to 0$ and $ta_t\to\infty$,
	the family of random variables $\left\{\frac{(a_tt)^{\alpha(\nu)}S(L_\nu(t))}{t}:t>0\right\}$ 
	satisfies the LDP with speed $1/a_t$ and good rate function $I_{\mathrm{MD}}(\cdot;\nabla\kappa_S(0))$ 
	defined by
	$$I_{\mathrm{MD}}(x;\nabla\kappa_S(0)):=
	\sup_{\theta\in\mathbb{R}^h}\{\langle\theta,x\rangle-\Lambda_{\nu,\nabla\kappa_S(0)}(\theta)\},$$
	where
	$$\Lambda_{\nu,\nabla\kappa_S(0)}(\theta):=\left\{\begin{array}{ll}
		\left(\frac{\langle\theta,Q\theta\rangle}{2}\right)^{1/\nu}&\ \mbox{if}\ \nabla\kappa_S(0)=0\\
		(\langle\theta,\nabla\kappa_S(0)\rangle)^{1/\nu}1_{\langle\theta,\nabla\kappa_S(0)\rangle\geq 0}&
		\ \mbox{if}\ \nabla\kappa_S(0)\neq 0.
	\end{array}\right.$$
\end{proposition}
\begin{proof}
	We consider an application of the G\"artner Ellis Theorem and we distinguish the cases 
	$\nabla\kappa_S(0)=0$ and $\nabla\kappa_S(0)\neq 0$ (in both cases the function 
	$\Lambda_{\nu,\nabla\kappa_S(0)}$ is finite and differentiable, and therefore we can apply that theorem).
	So we have to show that
	$$\lim_{t\to\infty}
	\frac{1}{1/a_t}\log\mathbb{E}\left[e^{\frac{1}{a_t}\langle\theta,\frac{(a_tt)^{\alpha(\nu)}S(L_\nu(t))}{t}\rangle}\right]
	=\Lambda_{\nu,\nabla\kappa_S(0)}(\theta)\ \mbox{for all}\ \theta\in\mathbb{R}^h,$$
	or equivalently (this can be checked with some standard computations)
	$$\lim_{t\to\infty}a_t\log E_\nu\left(\kappa_S\left(\frac{\theta}{(a_tt)^{1-\alpha(\nu)}}\right)t^\nu\right)
	=\Lambda_{\nu,\nabla\kappa_S(0)}(\theta)\ \mbox{for all}\ \theta\in\mathbb{R}^h.$$
	
	If $\nabla\kappa_S(0)=0$ we have
	\begin{multline*}
		a_t\log E_\nu\left(\kappa_S\left(\frac{\theta}{(a_tt)^{1-\alpha(\nu)}}\right)t^\nu\right)
		=a_t\log E_\nu\left(\left(\frac{\langle\theta,Q\theta\rangle}{2(a_tt)^\nu}
		+o\left(\frac{1}{(a_tt)^\nu}\right)\right)t^\nu\right)\\
		=a_t\log E_\nu\left(\frac{1}{a_t^\nu}\left(\frac{\langle\theta,Q\theta\rangle}{2}
		+(a_tt)^\nu o\left(\frac{1}{(a_tt)^\nu}\right)\right)\right)
		\to\left(\frac{\langle\theta,Q\theta\rangle}{2}\right)^{1/\nu}\ \mbox{for all}\ \theta\in\mathbb{R}^h.
	\end{multline*}
	
	If $\nabla\kappa_S(0)\neq 0$ we have
	\begin{multline*}
		a_t\log E_\nu\left(\kappa_S\left(\frac{\theta}{(a_tt)^{1-\alpha(\nu)}}\right)t^\nu\right)
		=a_t\log E_\nu\left(\left(\frac{\langle\theta,\nabla\kappa_S(0)\rangle}{(a_tt)^\nu}
		+o\left(\frac{1}{(a_tt)^\nu}\right)\right)t^\nu\right)\\
		=a_t\log E_\nu\left(\frac{1}{a_t^\nu}\left(\langle\theta,\nabla\kappa_S(0)\rangle
		+(a_tt)^\nu o\left(\frac{1}{(a_tt)^\nu}\right)\right)\right)
		\to(\langle\theta,\nabla\kappa_S(0)\rangle)^{1/\nu}1_{\langle\theta,\nabla\kappa_S(0)\rangle\geq 0}
		\ \mbox{for all}\ \theta\in\mathbb{R}^h.
	\end{multline*}
\end{proof}

\begin{remark}\label{rem:fills-the-gap-IMM}
	The class of LDPs in Proposition \ref{prop:ncMD-IMM} fills the gap between two asymptotic regimes.
	\begin{itemize}
		\item If $a_t=1$ (thus $a_t\to 0$ fails), then we have the weak convergence in Proposition \ref{prop:weak-convergence-IMM}
		(actually a weak convergence result for $\kappa_S(0)=0$, and another one for $\kappa_S(0)\neq 0$).
		\item If $a_t=\frac{1}{t}$ (thus $ta_t\to\infty$ fails), then we have the case of the reference LDP in Proposition 
		\ref{prop:reference-LDP-IMM}.
	\end{itemize}
\end{remark}

\begin{remark}\label{rem:literature-IMM}
	Here we collect some remarks concerning the results in this section and other results in the literature.\\
	(i) One could consider a generalization of the reference LDP in Proposition \ref{prop:reference-LDP-IMM} by considering an 
	independent random time-change having the same behavior of the inverse stable subordinator. This aspect was already
	discussed in \cite{IulianoMacciMeoli} (more precisely see Proposition 5.1, together with Condition 5.1, in that reference).\\
	(ii) We can adapt the content of Remark 3.3 in \cite{IulianoMacciMeoli} for Proposition \ref{prop:ncMD-IMM} in this paper.
	Indeed, if $Q$ is the null matrix and $\nabla\kappa_S(0)$ is the null vector, the process $\{S(L_\nu(t)):t\geq 0\}$ in 
	Propositions \ref{prop:reference-LDP-IMM}, \ref{prop:weak-convergence-IMM}, \ref{prop:ncMD-IMM} is identically equal
	to the null vector (because $S(t)=0\in\mathbb{R}^h$ for all $t\geq 0$) and the weak convergence in Proposition 
	\ref{prop:weak-convergence-IMM} (for $\nabla\kappa_S(0)=0$) would be towards a constant random variable equal to the null
	vector $0\in\mathbb{R}^h$.\\
	(iii) We remark that, in this section, the independence of the marginal univariate Lévy processes 
	$\{S_1(t):t\geq 0\},\ldots,\{S_h(t):t\geq 0\}$ is not required. On the contrary the independence of those processes was 
	required for the results in \cite{GuptaMacci} (see Conditions 1.4 and 1.5, together with Assumption 1.1, in that reference).
\end{remark}

In general we cannot provide an explicit expression of the rate function $I_{\mathrm{MD}}(\cdot;\nabla\kappa_S(0))$ in
Proposition \ref{prop:ncMD-IMM}. However, as shown in Proposition 3.3 in \cite{IulianoMacciMeoli}, this is possible for $h=1$.
More precisely, if we use the notation $\kappa_S^\prime(0)$ instead of $\nabla\kappa_S(0)$, and $q>0$ in place of the matrix $Q$,
we have: if $\kappa_S^\prime(0)=0$,
$$I_{\mathrm{MD}}(x;\kappa_S^\prime(0))=
((\nu/2)^{\nu/(2-\nu)}-(\nu/2)^{2/(2-\nu)})\left(\frac{2x^2}{q}\right)^{1/(2-\nu)};$$
if $\kappa_S^\prime(0)\neq 0$,
\begin{equation}\label{eq:auxiliary-h=1}
	I_{\mathrm{MD}}(x;\kappa_S^\prime(0))=\left\{\begin{array}{ll}
		(\nu^{\nu/(1-\nu)}-\nu^{1/(1-\nu)})\left(\frac{x}{\kappa_S^\prime(0)}\right)^{1/(1-\nu)}&
		\ \mbox{if}\ \frac{x}{\kappa_S^\prime(0)}\geq 0\\
		\infty&\ \mbox{if}\ \frac{x}{\kappa_S^\prime(0)}<0
	\end{array}\right.=:H_\nu(x;\kappa_S^\prime(0)).
\end{equation}

Here we can only provide explicit expressions of $I_{\mathrm{MD}}(x;\nabla\kappa_S(0))$ for some values of 
$x=(x_1,\ldots,x_h)$ when $\nabla\kappa_S(0)\neq 0$. In particular we refer to the function 
$H_\nu(\cdot;\kappa_S^\prime(0))$ (see eq.~\eqref{eq:auxiliary-h=1}) concerning the case $h=1$.

\begin{proposition}\label{prop:explicit-cases-IMM}
	Assume that $\nabla\kappa_S(0)=(m_1,\ldots,m_h)\neq 0$. Then we have the following cases.\\
	$(i)$ If there exists $i\in\{1,\ldots,h\}$ such that $x_im_i<0$, then
	$I_{\mathrm{MD}}(x;\nabla\kappa_S(0))=\infty$.\\
	$(ii)$ If there exists $i\in\{1,\ldots,h\}$ such that $m_i=0$ and $x_i\neq 0$, then
	$I_{\mathrm{MD}}(x;\nabla\kappa_S(0))=\infty$.\\
	$(iii)$ If there exists $c(x)\geq 0$ such that $x_i=c(x)m_i$ for every 
	$i\in\{1,\ldots,h\}$, then $I_{\mathrm{MD}}(x;\nabla\kappa_S(0))=H_\nu(c(x);1)$
	or, equivalently, $I_{\mathrm{MD}}(x;\nabla\kappa_S(0))=H_\nu(x_i;m_i)$ for every 
	$i\in\{1,\ldots,h\}$.
\end{proposition}
\begin{proof}
	We prove the three statements separately.\\
	If there exists $i\in\{1,\ldots,h\}$ such that $x_im_i<0$, we prove $(i)$ noting that
	$$I_{\mathrm{MD}}(x;\nabla\kappa_S(0))\geq
	\sup_{\theta_i\in\mathbb{R}}\{\theta_ix_i-(\theta_im_i)^{1/\nu}1_{\theta_im_i\geq 0}\}.$$
	Indeed, the supremum at the right hand side is equal to $+\infty$ by taking the limit as $\theta_i\to+\infty$
	(when $x_i>0$ and $m_i<0$), or by taking the limit as $\theta_i\to-\infty$ (when $x_i<0$ and $m_i>0$).\\
	If there exists $i\in\{1,\ldots,h\}$ such that $m_i=0$ and $x_i\neq 0$, then we prove $(ii)$ noting that
	$$I_{\mathrm{MD}}(x;\nabla\kappa_S(0))\geq
	\sup_{\theta_i\in\mathbb{R}}\{\theta_ix_i-(\theta_im_i)^{1/\nu}1_{\theta_im_i\geq 0}\}
	=\sup_{\theta_i\in\mathbb{R}}\theta_ix_i=+\infty.$$
	If there exists $c(x)\geq 0$ such that $x_i=c(x)m_i$ for every $i\in\{1,\ldots,h\}$, then we prove $(iii)$ 
	noting that
	$$I_{\mathrm{MD}}(x;\nabla\kappa_S(0))=
	\sup_{\theta\in\mathbb{R}^h}\{c(x)\langle\theta,\nabla\kappa_S(0)\rangle-
	(\langle\theta,\nabla\kappa_S(0)\rangle)^{1/\nu}1_{\langle\theta,\nabla\kappa_S(0)\rangle\geq 0}\}
	=\sup_{\eta\in\mathbb{R}}\{c(x)\eta-\eta^{1/\nu}1_{\eta\geq 0}\}.$$
	Indeed, one can check that the supremum is attained at $\eta=(\nu c(x))^{\nu/(1-\nu)}$, and we get
	$$I_{\mathrm{MD}}(x;\nabla\kappa_S(0))=(\nu c(x))^{\nu/(1-\nu)}c(x)-(\nu c(x))^{1/(1-\nu)}=
	(\nu^{\nu/(1-\nu)}-\nu^{1/(1-\nu)})(c(x))^{1/(1-\nu)},$$
	where the last expression coincides with $H_\nu(c(x);1)$ or, equivalently, with $H_\nu(x_i;m_i)$ for every 
	$i\in\{1,\ldots,h\}$.
\end{proof}

\section{Multivariate Lévy processes time-changed with a subordinator}\label{sec:Levy}
Let $\{S(t):t\geq 0\}$ be a $\mathbb{R}^h$-valued Lévy process as in Condition \ref{cond:IMM-Levy}.
Moreover, let $\{V(t):t\geq 0\}$ be a subordinator such that $\kappa_V(\eta):=\log\mathbb{E}\left[e^{\eta V(1)}\right]$ 
is finite when $\eta$ belongs to a neighborhood of the origin $0\in\mathbb{R}$. It is well-known that
$$\mathbb{E}[V(1)]=\kappa_V^\prime(0).$$
We assume that $\{V(t):t\geq 0\}$ is not trivial, i.e. $\kappa_V^\prime(0)>0$. Indeed, if we had $\kappa_V^\prime(0)=0$,
then we would have $V(t)=0$ and $S(V(t))=0$ for every $t\geq 0$, and this is not interesting.
In what follows we always assume that $\{S(t):t\geq 0\}$ and $\{V(t):t\geq 0\}$ are independent.

It is well-known that $\kappa_S$ and $\kappa_V$ are lower semicontinuous functions (see, e.g.,~Exercise 2.2.22 in 
\cite{DemboZeitouni}). Then the function $\kappa_V\circ\kappa_S$ is also lower semicontinuous since it is the
logarithm of a moment generating function. Indeed,
$$\kappa_V\circ\kappa_S(\theta)=\log\mathbb{E}\left[e^{\langle\theta,S(V(1))\rangle}\right].$$

We have the following results.

\begin{proposition}[Reference LDP]\label{prop:reference-LDP-Levy}
    Assume that $\kappa_V\circ\kappa_S$ is an essentially smooth function. Then 
	$\left\{\frac{S(V(t))}{t}:t>0\right\}$ satisfies the LDP with speed $t$ and good rate function $I_{\mathrm{LD}}$ 
	defined by
	$$I_{\mathrm{LD}}(x):=\sup_{\theta\in\mathbb{R}^h}\{\langle\theta,x\rangle-\kappa_V(\kappa_S(\theta))\}.$$
\end{proposition}
\begin{proof}
	We consider a straightforward application of the G\"artner Ellis Theorem. Indeed, we have (note that 
	$\kappa_V(\kappa_S(\theta))=+\infty$ when $\kappa_S(\theta)=+\infty$)
	$$\frac{1}{t}\log\mathbb{E}\left[e^{t\langle\theta,S(V(t))/t\rangle}\right]
	=\frac{1}{t}\log\mathbb{E}\left[e^{V(t)\kappa_S(\theta)}\right]=\kappa_V(\kappa_S(\theta))
	\quad(\mbox{for all}\ \theta\in\mathbb{R}^h).$$
\end{proof}

\begin{proposition}[Weak convergence]\label{prop:weak-convergence-Levy}
	$\left\{S\left(\frac{V(t)}{t}\right):t>0\right\}$ converges weakly to $S(\kappa_V^\prime(0))$ as $t\to\infty$.
	More precisely, we mean
	$$\lim_{t\to\infty}\mathbb{E}\left[e^{\langle\theta,S\left(\frac{V(t)}{t}\right)\rangle}\right]=
    \mathbb{E}\left[e^{\langle\theta,S(\kappa_V^\prime(0))\rangle}\right]=e^{\kappa_V^\prime(0)\kappa_S(\theta)}
	\quad(\mbox{for all}\ \theta\in\mathbb{R}^h).$$
\end{proposition}
\begin{proof}
	The desired limit trivially holds if $\kappa_S(\theta)=+\infty$. On the contrary, if 
	$\kappa_S(\theta)<+\infty$, for $t$ large enough $\frac{\kappa_S(\theta)}{t}$ is close to zero and 
	$\kappa_V\left(\frac{\kappa_S(\theta)}{t}\right)$ is also finite. Thus
	$$\mathbb{E}\left[e^{\langle\theta,S\left(\frac{V(t)}{t}\right)\rangle}\right]
	=\mathbb{E}\left[e^{\frac{V(t)}{t}\kappa_S(\theta)}\right]
	=\exp\left(t\kappa_V\left(\frac{\kappa_S(\theta)}{t}\right)\right)
	\to\exp\left(\kappa_V^\prime(0)\kappa_S(\theta)\right)
	\quad(\mbox{for all}\ \theta\in\mathbb{R}^h).$$
	Indeed, we take into account that $\kappa_V(0)=0$ and therefore
	$$t\kappa_V\left(\frac{\kappa_S(\theta)}{t}\right)=
	\frac{\kappa_V\left(\frac{\kappa_S(\theta)}{t}\right)-\kappa_V(0)}{\kappa_S(\theta)/t}\kappa_S(\theta)
	\to\kappa_V^\prime(0)\kappa_S(\theta).$$
\end{proof}

\begin{proposition}[Noncentral moderate deviations]\label{prop:ncMD-Levy}
	Assume that $\kappa_S$ is an essentially smooth function. Then, for every $a_t>0$ such that
	$a_t\to 0$ and $ta_t\to\infty$, $\left\{a_tS\left(\frac{V(t)}{ta_t}\right):t>0\right\}$
	satisfies the LDP with speed $1/a_t$ and good rate function $I_{\mathrm{MD}}$ defined by
	$$I_{\mathrm{MD}}(x):=\sup_{\theta\in\mathbb{R}^h}\left\{\langle\theta,x\rangle-\kappa_V^\prime(0)\kappa_S(\theta)\right\}.$$
\end{proposition}
\begin{proof}
	We consider a straightforward application of the G\"artner Ellis Theorem. Indeed, we have 
	(note that $\kappa_V\left(\frac{\kappa_S(\theta)}{ta_t}\right)=+\infty$ when $\kappa_S(\theta)=+\infty$; 
	moreover, when $\kappa_S(\theta)<+\infty$, $\frac{\kappa_S(\theta)}{ta_t}$ is close to zero 
	for $t$ large enough, and therefore $\kappa_V\left(\frac{\kappa_S(\theta)}{ta_t}\right)$ is finite)	
	\begin{multline*}
		\frac{1}{1/a_t}\log\mathbb{E}\left[e^{\langle\frac{\theta}{a_t},a_tS\left(\frac{V(t)}{ta_t}\right)\rangle}\right]
		=a_t\log\mathbb{E}\left[e^{\langle\theta,S\left(\frac{V(t)}{ta_t}\right)\rangle}\right]\\
		=a_t\log\mathbb{E}\left[e^{\frac{V(t)}{ta_t}\kappa_S(\theta)}\right]
		=a_tt\kappa_V\left(\frac{\kappa_S(\theta)}{ta_t}\right)\to\kappa_V^\prime(0)\kappa_S(\theta)
		\quad(\mbox{for all}\ \theta\in\mathbb{R}^h),
	\end{multline*}
    where we take into account that $\kappa_V(0)=0$ and $ta_t\to\infty$. Moreover, $\kappa_V^\prime(0)\kappa_S(\cdot)$
    is an essentially smooth function (by the essential smoothness of $\kappa_S$). 
\end{proof}

\begin{remark}\label{rem:fills-the-gap-Levy}
	The class of LDPs in Proposition \ref{prop:ncMD-Levy} fills the gap between two asymptotic regimes.
	\begin{itemize}
		\item If $a_t=1$ (thus $a_t\to 0$ fails), then we have the weak convergence in Proposition \ref{prop:weak-convergence-Levy}
		(note that we have a non-central result if $\{S(t):t\geq 0\}$ is not a Brownian motion, possibly with drift).
		\item If $a_t=\frac{1}{t}$ (thus $ta_t\to\infty$ fails), then we have the case of the reference LDP in Proposition 
		\ref{prop:reference-LDP-Levy}.
	\end{itemize}
\end{remark}

\begin{remark}\label{rem:inequalities-between-rfs-Levy}
	By taking into account $\kappa_V(0)=0$ and the convexity of $\kappa_V$, we have 
	$\kappa_V^\prime(0)\kappa_S(\theta)\leq\kappa_V(\kappa_S(\theta))$ for all $\theta\in\mathbb{R}^h$. So we have
	\begin{equation}\label{eq:inequality-Levy}
		I_{\mathrm{MD}}(x)\geq I_{\mathrm{LD}}(x)\quad(\mbox{for all}\ x\in\mathbb{R}^h).
	\end{equation}
    Note that both rate functions $I_{\mathrm{LD}}$ and $I_{\mathrm{MD}}$ uniquely vanish at 
    $x=\kappa_V^\prime(0)\nabla\kappa_S(0)$.
\end{remark}

\section{On some examples with a compound Poisson weak limits}\label{sec:Poisson}
Let $\{X_n(p):n\geq 1,p\in[0,1]\}$ be a sequence of i.i.d. $\mathbb{R}^h$-valued random variables such that
$$\mathbb{E}\left[e^{\langle\theta,X_n(p)\rangle}\right]=1-p+pG(\theta)\quad(\mbox{for all}\ \theta\in\mathbb{R}^h),$$
where $p=P(X_n(p)=0)$ and $G(\theta)=\mathbb{E}[e^{\langle\theta,X_n(p)\rangle}|X_n(p)\neq 0]$. So we are assuming that
the conditional distribution of $X_n(p)$ given $X_n(p)\neq 0$ does not depend on $p$ (and the moment generating function
of this conditional distribution is $G$). We also assume that $G$ is finite in a neighborhood of the origin 
$0\in\mathbb{R}^h$.

In this section we consider sequences $\{p_n:n\geq 1\}\subset[0,1]$ in place of $p$. In Propositions 
\ref{prop:weak-convergence-Poisson} and \ref{prop:ncMD-Poisson} we have $p_n\to 0$ as a consequence of the 
hypotheses (i.e. $np_n\to\lambda\in(0,\infty)$, and $na_np_n\to\lambda\in(0,\infty)$ when $na_n\to\infty$, respectively).
In Proposition \ref{prop:reference-LDP-Poisson} we have $p_n\to p\in[0,1]$ and, in order to avoid trivialities, we think
to have $p_n\to p\in(0,1)$.

We have the following results.

\begin{proposition}[Reference LDP]\label{prop:reference-LDP-Poisson}
	Let $\{p_n:n\geq 1\}$ be a sequence such that $p_n\to p\in(0,1)$. Moreover, if $\{p_n:n\geq 1\}$ is not a constant
	sequence, assume that $G$ is an essentially smooth function. Then $\left\{\frac{X_1(p_n)+\cdots+X_n(p_n)}{n}:n\geq 1\right\}$
	satisfies the LDP with speed $n$ and good rate function $I_{\mathrm{LD}}$ defined by
	$$I_{\mathrm{LD}}(x):=\sup_{\theta\in\mathbb{R}^h}\{\langle\theta,x\rangle-\log(1+p(G(\theta)-1))\}.$$
\end{proposition}
\begin{proof}
	If the sequence $\{p_n:n\geq 1\}$ is constant, then we consider a straightforward application of the Cramér Theorem.
	On the contrary, if $\{p_n:n\geq 1\}$ is not a constant sequence, then we consider a straightforward application of 
	the G\"artner Ellis Theorem. Indeed, we have (note that, when $G(\theta)=+\infty$, the functions of $\theta$ that 
	appear below are equal to infinity)
	$$\frac{1}{n}\log\mathbb{E}\left[e^{\langle\theta,X_1(p_n)+\cdots+X_n(p_n)\rangle}\right]
	=\log(1+p_n(G(\theta)-1))\to\log(1+p(G(\theta)-1))
	\quad(\mbox{for all}\ \theta\in\mathbb{R}^h)$$
    and, moreover, $\log(1+p(G(\cdot)-1))$ is an essentially smooth function (by the essential smoothness of $G$).
\end{proof}

\begin{proposition}[Weak Convergence]\label{prop:weak-convergence-Poisson}
	Let $\{p_n:n\geq 1\}$ be a sequence such that $np_n\to\lambda\in(0,\infty)$. Then
	$\left\{X_1(p_n)+\cdots+X_n(p_n):n\geq 1\right\}$ converges weakly to a compound Poisson distributed random variable $Y$
	such that the intensity of the jumps is $\lambda$ and the moment generating function of the jumps is $G(\theta)$. More
	precisely, we mean that
	$$\lim_{n\to\infty}\mathbb{E}\left[e^{\langle\theta,X_1(p_n)+\cdots+X_n(p_n)\rangle}\right]=
	\mathbb{E}\left[e^{\langle\theta,Y\rangle}\right]=\exp(\lambda(G(\theta)-1))\quad(\mbox{for all}\ \theta\in\mathbb{R}^h).$$
\end{proposition}
\begin{proof}
	The desired limit trivially holds if $G(\theta)=+\infty$. On the contrary, if $G(\theta)<+\infty$, for $n$ large enough 
	$(1+p_n(G(\theta)-1))$ is also finite. Thus
	\begin{multline*}
		\mathbb{E}\left[e^{\langle\theta,X_1(p_n)+\cdots+X_n(p_n)\rangle}\right]=(1+p_n(G(\theta)-1))^n\\
		=\exp\left(np_n\frac{\log(1+p_n(G(\theta)-1))}{p_n}\right)\to\exp(\lambda(G(\theta)-1))
		\quad(\mbox{for all}\ \theta\in\mathbb{R}^h).
	\end{multline*}
\end{proof}

\begin{proposition}[Noncentral moderate deviations]\label{prop:ncMD-Poisson}
	Assume that $G$ is an essentially smooth function. Then, for every $a_n>0$ such that
	$a_n\to 0$ and $na_n\to\infty$, and also $na_np_n\to\lambda\in(0,\infty)$, 
	$\left\{\frac{X_1(p_n)+\cdots+X_n(p_n)}{1/a_n}:n\geq 1\right\}$ satisfies the LDP with speed $1/a_n$ and good rate function 
	$I_{\mathrm{MD}}$ defined by
	$$I_{\mathrm{MD}}(x):=\sup_{\theta\in\mathbb{R}^h}\left\{\langle\theta,x\rangle-
	\lambda(G(\theta)-1)\right\}.$$
\end{proposition}
\begin{proof}
	We consider a straightforward application of the G\"artner Ellis Theorem. Indeed, we have (note that, 
	when $G(\theta)=+\infty$, the functions of $\theta$ that appear below are equal to infinity)
	\begin{multline*}
		\frac{1}{1/a_n}\log\mathbb{E}\left[e^{\langle\theta,X_1(p_n)+\cdots+X_n(p_n)\rangle}\right]=na_n\log(1+p_n(G(\theta)-1))\\
		=na_np_n\frac{\log(1+p_n(G(\theta)-1))}{p_n}\to\lambda(G(\theta)-1)
		\quad(\mbox{for all}\ \theta\in\mathbb{R}^h).
	\end{multline*}
    and, moreover, $\lambda(G(\cdot)-1)$ is an essentially smooth function (by the essential smoothness of $G$).
\end{proof}

We have the following remarks.

\begin{remark}\label{rem:fills-the-gap-Poisson}
	The class of LDPs in Proposition \ref{prop:ncMD-Poisson} fills the gap between two asymptotic regimes.
	\begin{itemize}
		\item If $a_n=1$ (thus $a_n\to 0$ fails), then we have the weak convergence in Proposition 
		\ref{prop:weak-convergence-Poisson}.
		\item If $a_n=\frac{1}{n}$ (thus $na_n\to\infty$ fails) and $\lambda\in(0,1)$, then we have the case of 
		the reference LDP in Proposition \ref{prop:reference-LDP-Poisson} with $\lambda=p$.
	\end{itemize}
\end{remark}

\begin{remark}\label{rem:inequalities-between-rfs-Poisson}
	It is easy to check that we have $\log(1+p(G(\theta)-1))\leq p(G(\theta)-1)$ for all $\theta\in\mathbb{R}^h$. So, if we 
	take $\lambda=p$ (see the second item in Remark \ref{rem:fills-the-gap-Poisson}), then we have
	\begin{equation}\label{eq:inequality-Poisson}
		I_{\mathrm{LD}}(x)\geq I_{\mathrm{MD}}(x)\quad(\mbox{for all}\ x\in\mathbb{R}^h).
	\end{equation}
	Note that, since $\lambda=p$, both rate functions $I_{\mathrm{LD}}$ and $I_{\mathrm{MD}}$ uniquely vanish at the same point,
	i.e. $x=\lambda\nabla G(0)$ or $x=p\nabla G(0)$.
\end{remark}

\begin{remark}\label{rem:on-the-well-known-case-Poisson}
	The inequality \eqref{eq:inequality-Poisson} can be checked for the case $h=1$ and $G(\theta)=e^\theta$ (the well-known case 
	of the weak convergence of the Binomial distributions to the Poisson distribution). In this case we have closed form expressions
	of the rate functions, i.e.
	$$I_{\mathrm{LD}}(x)=\left\{\begin{array}{ll}
		x\log\frac{x}{p}+(1-x)\log\frac{1-x}{1-p}&\ \mbox{for}\ x\in[0,1]\\
		\infty&\ \mbox{otherwise}
	\end{array}\right.$$
	and, since $\lambda=p$ (see Remarks \ref{rem:fills-the-gap-Poisson} and \ref{rem:inequalities-between-rfs-Poisson}),
	$$I_{\mathrm{MD}}(x)=\left\{\begin{array}{ll}
		x\log\frac{x}{p}-x+p&\ \mbox{for}\ x\geq 0\\
		\infty&\ \mbox{otherwise}.
	\end{array}\right.$$
	Then the inequality \eqref{eq:inequality-Poisson} is trivial if $x\notin [0,1]$. Otherwise, for $x\in [0,1]$, it is 
	easy to check that the function
	$$x\mapsto I_{\mathrm{LD}}(x)-I_{\mathrm{MD}}(x)=(1-x)\log\frac{1-x}{1-p}+x-p$$
	attains its minimum at $x=p$ and it is equal to zero because $I_{\mathrm{LD}}(p)-I_{\mathrm{MD}}(p)=0-0=0$.
\end{remark}

\section{On some examples obtained by applying continuous mappings}\label{sec:contraction-continuity}
A well-known (central) moderate deviation result is Theorem 3.7.1 in \cite{DemboZeitouni} which concerns sums of i.i.d. 
$\mathbb{R}^h$ valued random variables. More precisely, under suitable hypotheses, this result fills the gap between the convergence 
to a constant (provided by the law of large numbers on $\mathbb{R}^h$) and a weak convergence to a Normal distribution (provided by 
the Central Limit Theorem on $\mathbb{R}^h$). Moreover, the reference LDP for the convergence to a constant is provided by the Cramér 
Theorem (see, e.g.,~Theorem 2.2.3 and Theorem 2.2.30 in \cite{DemboZeitouni} with the subsequent Remark (a)).

The aim of this section is to present a family of other (possibly noncentral) moderate deviation results derived by Theorem 3.7.1 in 
\cite{DemboZeitouni} and by applying continuous transformations $U:D\subset\mathbb{R}^h\to\mathbb{R}^k$ for some subset $D$.
This can be done by applying the Contraction Principle and the well-known Continuity Theorem for the weak convergence. In particular
this trick can be considered to obtain other moderate deviation results. In other words it is possible to present
other (possibly noncentral) moderate deviation results by applying a continuous function $U$ to the random variables that appear 
in a (possibly noncentral) moderate deviation result different from Theorem 3.7.1 in \cite{DemboZeitouni}.

In view of what follows we introduce some notation. Let $\{X_n:n\geq 1\}$ be a sequence of $\mathbb{R}^h$-valued random variables 
such that
$$\kappa_X(\theta):=\log\mathbb{E}\left[e^{\langle\theta,X_1\rangle}\right]$$
is finite when $\theta$ belongs to a neighborhood of the origin $0\in\mathbb{R}^h$. We recall that
$$\mathbb{E}[X_1]=\nabla\kappa_X(0)\quad\mbox{and}\quad (\mathrm{Cov}(X_1^{(i)},X_1^{(j)}))_{i,j}=H\kappa_X(0),$$
where $H\kappa_X(0)$ is the Hessian matrix of $\kappa_X$ at the origin $\theta=0$.

We have the following results.

\begin{proposition}[Reference LDP]\label{prop:reference-LDP-continuity}
	The sequence $\left\{U\left(\frac{X_1+\cdots+X_n}{n}\right):n\geq 1\right\}$ satisfies 
	the LDP with speed $n$ and good rate function $I_{\mathrm{LD},U}$ defined by
	$$I_{\mathrm{LD},U}(y):=\inf_{x\in\mathbb{R}^h}\{\kappa_X^*(x):U(x)=y\},$$
	where
	$$\kappa_X^*(x):=\sup_{\theta\in\mathbb{R}^h}\{\langle\theta,x\rangle-\kappa_X(\theta)\}.$$
\end{proposition}
\begin{proof}
	It is an immediate consequence of a straightforward application of the Cramér Theorem combined with the Contraction Principle.
\end{proof}

\begin{proposition}[Weak Convergence]\label{prop:weak-convergence-continuity}
	The sequence $\left\{U\left(\frac{X_1+\cdots+X_n-n\nabla\kappa_X(0)}{\sqrt{n}}\right):n\geq 1\right\}$
	converges weakly to $U(Z)$, where $Z$ is a centered Normal distribution with covariance matrix $H\kappa_X(0)$.
\end{proposition}
\begin{proof}
	It is an immediate consequence of a straightforward application of the Central Limit Theorem combined with the Continuity Theorem.
\end{proof}

\begin{proposition}[Noncentral moderate deviations]\label{prop:ncMD-continuity}
	For every $a_n>0$ such that $a_n\to 0$ and $na_n\to\infty$, 
	$\left\{U\left(\frac{X_1+\cdots+X_n-n\nabla\kappa_X(0)}{\sqrt{n/a_n}}\right):n\geq 1\right\}$
	satisfies the LDP with speed $1/a_n$ and good rate function $I_{\mathrm{MD},U}$ defined by
	$$I_{\mathrm{MD},U}(y):=\inf_{x\in\mathbb{R}^h}\{\widetilde{\kappa}_X^*(x):U(x)=y\},$$
	where
	$$\widetilde{\kappa}_X^*(x):=\sup_{\theta\in\mathbb{R}^h}\left\{\langle\theta,x\rangle
	-\frac{\langle\theta,H\kappa_X(0)\theta\rangle}{2}\right\}.$$
\end{proposition}
\begin{proof}
	For every choice of $\{a_n:n\geq 1\}$, the desired LDP is an immediate consequence of a straightforward application of Theorem 
	3.7.1 in \cite{DemboZeitouni} combined with the Contraction Principle.
\end{proof}

We have the following remarks.

\begin{remark}\label{rem:fills-the-gap-contraction-continuity}
	The class of LDPs in Proposition \ref{prop:ncMD-continuity} fills the gap between two asymptotic regimes.
	\begin{itemize}
		\item If $a_n=1$ (thus $a_n\to 0$ fails), then we have the weak convergence in Proposition 
		\ref{prop:weak-convergence-continuity} (note that we have a non-central result if $U(Z)$ in Proposition 
		\ref{prop:weak-convergence-continuity} is not Gaussian distributed).
		\item If $a_n=\frac{1}{n}$ (thus $na_n\to\infty$ fails), then we have the case of the reference LDP in Proposition 
		\ref{prop:reference-LDP-continuity} with a possible shift. More precisely, we mean the LDP of $\left\{U\left(\frac{X_1+\cdots+X_n}{n}-\nabla\kappa_X(0)\right):n\geq 1\right\}$ with speed $n$ and good rate function $\widetilde{I}_{\mathrm{LD},U}$ defined by
		$$\widetilde{I}_{\mathrm{LD},U}(y):=\inf_{x\in\mathbb{R}^h}\{\kappa_X^*(x+\nabla\kappa_X(0)):U(x)=y\}.$$
	\end{itemize}
\end{remark}

\begin{remark}\label{rem:various-continuity}
	The rate function $I_{\mathrm{LD},U}$ in Proposition \ref{prop:reference-LDP-continuity} uniquely vanishes at 
	$y=U(\nabla\kappa_X(0))$. The rate function $\widetilde{I}_{\mathrm{MD},U}$ in Proposition \ref{prop:ncMD-continuity} uniquely
	vanishes at $y=U(0)$. Moreover, it is well-known that, if $H\kappa_X(0)$ is invertible, in Proposition \ref{prop:ncMD-continuity} we have
	$$\widetilde{\kappa}_X^*(x)=\frac{\langle x,(H\kappa_X(0))^{-1}x\rangle}{2}$$
	(where $\langle\cdot,\cdot\rangle$ is the inner product in $\mathbb{R}^h$).
\end{remark}

In what follows we present some examples. In particular we emphasize the case in which the function $U$ is vector-valued, i.e. $k\geq 2$.
However in the examples above the case $k=1$ is also allowed. In particular we discuss the possibility to get explicit formulas for
$I_{\mathrm{LD},U}$ and $I_{\mathrm{MD},U}$.

\begin{example}\label{ex1}
	Here we consider the function $U_1:\mathbb{R}^h\to\mathbb{R}^{h+1}$ defined by
	$$U_1(x_1,\ldots,x_h):=\left(\frac{e^{x_1}}{1+\sum_{j=1}^he^{x_j}},\ldots,\frac{e^{x_h}}{1+\sum_{j=1}^he^{x_j}},
	\frac{1}{1+\sum_{j=1}^he^{x_j}}\right).$$
	In this case $U_1(Z)$ has \emph{logistic Normal distribution}, and the range of $U_1$ is the simplex
	$$\left\{(y_1,\ldots,y_{h+1})\in\mathbb{R}^{h+1}:y_1,\ldots,y_{h+1}\geq 0,\sum_{j=1}^{h+1}y_j=1\right\}.$$
	Then we have
	$$I_{\mathrm{LD},U_1}(y_1,\ldots,y_{h+1}):=\kappa_X^*(\log(y_1/y_{h+1}),\ldots,\log(y_h/y_{h+1}))$$
	and
	$$I_{\mathrm{MD},U_1}(y_1,\ldots,y_{h+1}):=\widetilde{\kappa}_X^*(\log(y_1/y_{h+1}),\ldots,\log(y_h/y_{h+1})).$$
\end{example}

\begin{example}\label{ex2}
	Here we take $h\geq 2$ and
	$$H\kappa_X(0)=\left(\begin{array}{cc}
		\Psi&0\\
		0&1
	\end{array}\right),$$
	where $\Psi$ is a full range covariance matrix of order $h-1$ (thus $H\kappa_X(0)$ is also a full range covariance matrix). Then we 
	consider $U_2:\mathbb{R}^h\to\mathbb{R}^{h-1}$ defined by
	$$U_2(x_1,\ldots,x_h):=\left((1-\delta_1^2)^{1/2}x_1+\delta_1|x_h|,\ldots,(1-\delta_{h-1}^2)^{1/2}x_{h-1}+\delta_{h-1}|x_h|\right),$$
	for $\delta_1,\ldots,\delta_{h-1}\in(-1,1)$. In this case $U_2(Z)$ has a $(h-1)$-variate \emph{skew Normal distribution} (note that 
	$U_2(Z)$ is Gaussian distributed only if $\delta_1=\cdots=\delta_{h-1}=0$). Then we have
	$$I_{\mathrm{LD},U_2}(y_1,\ldots,y_{h-1}):=\inf_{x_h\in\mathbb{R}}\kappa_X^*\left(\frac{y_1-|x_h|\delta_1}{(1-\delta_1^2)^{1/2}},
	\ldots,\frac{y_{h-1}-|x_h|\delta_{h-1}}{(1-\delta_{h-1}^2)^{1/2}},x_h\right)$$
	and
	\begin{equation}\label{eq:MD-rf-skew}
		I_{\mathrm{MD},U_2}(y_1,\ldots,y_{h-1}):=\inf_{x_h\in\mathbb{R}}\widetilde{\kappa}_X^*
		\left(\frac{y_1-|x_h|\delta_1}{(1-\delta_1^2)^{1/2}},
		\ldots,\frac{y_{h-1}-|x_h|\delta_{h-1}}{(1-\delta_{h-1}^2)^{1/2}},x_h\right).
	\end{equation}
\end{example}

Actually, as we shall see, we can provide an explicit formula for the rate function $I_{\mathrm{MD},U_2}$ in \eqref{eq:MD-rf-skew}. 
We start noting that $H\kappa_X(0)$ is invertible and we have
$$(H\kappa_X(0))^{-1}=\left(\begin{array}{cc}
	\Psi^{-1}&0\\
	0&1
\end{array}\right).$$
Moreover, if we consider $a^{(y,\delta)},b^{(\delta)}\in\mathbb{R}^{h-1}$ defined by
$$a^{(y,\delta)}:=\left(\frac{y_j}{(1-\delta_j^2)^{1/2}}\right)_{j=1,\ldots,h-1}\ \mbox{and}\
b^{(\delta)}:=\left(\frac{\delta_j}{(1-\delta_j^2)^{1/2}}\right)_{j=1,\ldots,h-1},$$
then \eqref{eq:MD-rf-skew} yields 
$$I_{\mathrm{MD},U_2}(y_1,\ldots,y_{h-1}):=\inf_{x_h\in\mathbb{R}}\widetilde{\kappa}_X^*
\left(a^{(y,\delta)}-|x_h|b^{(\delta)},x_h\right),$$
and therefore (here $\langle\cdot,\cdot\rangle$ is the inner product in $\mathbb{R}^{h-1}$)
$$I_{\mathrm{MD},U_2}(y_1,\ldots,y_{h-1})=\frac{1}{2}\inf_{x_h\in\mathbb{R}}\left\{\langle a^{(y,\delta)},\Psi^{-1}a^{(y,\delta)}\rangle
+x_h^2\langle b^{(\delta)},\Psi^{-1}b^{(\delta)}\rangle+x_h^2-2|x_h|\langle a^{(y,\delta)},\Psi^{-1}b^{(\delta)}\rangle\right\}.$$
Note that we have the infimum of a parabola on the positive half-line and, moreover, the infimum over all the entire real line is 
attained at $x_h=\widehat{x}_h$, where
$$\widehat{x}_h:=\frac{\langle a^{(y,\delta)},\Psi^{-1}b^{(\delta)}\rangle}{\langle b^{(\delta)},\Psi^{-1}b^{(\delta)}\rangle+1}.$$
Therefore we have the following two cases:
\begin{itemize}
	\item $\langle a^{(y,\delta)},\Psi^{-1}b^{(\delta)}\rangle\leq 0$ (where the infimum over the positive half-line 
	is attained at $x_h=0$);
	\item $\langle a^{(y,\delta)},\Psi^{-1}b^{(\delta)}\rangle>0$ (where the infimum over the positive half-line is attained at 
	$x_h=\widehat{x}_h$).
\end{itemize}
In conclusion we have the following explicit expression for $I_{\mathrm{MD},U_2}$ in \eqref{eq:MD-rf-skew}:
$$I_{\mathrm{MD},U_2}(y_1,\ldots,y_{h-1})=\left\{\begin{array}{ll}
		\frac{1}{2}\langle a^{(y,\delta)},\Psi^{-1}a^{(y,\delta)}\rangle&\ \mbox{if}\ \langle a^{(y,\delta)},\Psi^{-1}b^{(\delta)}\rangle\leq 0\\
		\frac{1}{2}\left\{\langle a^{(y,\delta)},\Psi^{-1}a^{(y,\delta)}\rangle
		-\frac{(\langle a^{(y,\delta)},\Psi^{-1}b^{(\delta)}\rangle)^2}{\langle b^{(\delta)},\Psi^{-1}b^{(\delta)}\rangle
			+1}\right\}&\ \mbox{if}\ \langle a^{(y,\delta)},\Psi^{-1}b^{(\delta)}\rangle>0.
\end{array}\right.$$

\subsection*{Acknowledgements}
The authors thank two referees for a careful reading of the first version of this manuscript. Moreover, they also
thank Brunero Liseo for some discussions on the theory of skew distributions and for suggesting the reference
\cite{AzzaliniCapitanioJRSSB} in Example \ref{ex2}.

\end{document}